\documentclass[preprint,12pt]{elsarticle}

\usepackage{amsthm}
\usepackage{amssymb,amsmath}

\theoremstyle{plain}
\newtheorem{thm}{Theorem}
\newtheorem{lem}[thm]{Lemma}
\newtheorem{prop}[thm]{Proposition}
\newtheorem{cor}[thm]{Corollary}

\newtheorem{req}[thm]{Remark}
\newtheorem{nt}[thm]{Notation}

\theoremstyle{definition}

\newtheorem{prob}{Problem}

\newproof{pf}{Proof}
\newproof{pot}{Proof of Theorem \ref{main theorem}}

 \makeatletter
\@addtoreset{thm}{chapter}
\@addtoreset{thm}{section}
\@addtoreset{thm}{subsection}

\makeatother

\begin{document}

\begin{frontmatter}

\title{About the spectra of a real nonnegative matrix and its signings}

\author[B]{Kawtar Attas}
\ead{kawtar.attas@gmail.com}

\author[B]{Abderrahim Boussa\"{\i}ri\corref{cor1}}
\ead{aboussairi@hotmail.com}

\author[B]{Mohamed Zaidi}	
\ead{zaidi.fsac@gmail.com}

\cortext[cor1]{Corresponding author}

\address[B]{Facult\'e des Sciences A\"{\i}n Chock, D\'epartement de
Math\'ematiques et Informatique,

Km 8 route d'El Jadida,
BP 5366 Maarif, Casablanca, Maroc}

\begin{abstract}
For a real matrix $M$, we denote by $sp(M)$ the spectrum of $M$ and by
$\left \vert M\right \vert $ its absolute value, that is the matrix obtained
from $M$ by replacing each entry of $M$ by its absolute value. Let $A$ be a
nonnegative real matrix, we call a \emph{signing} of $A$ every real matrix $B$ such that $\left \vert
B\right \vert =A$. In this paper, we study
the set of all signings of $A$ such that $sp(B)=\alpha sp(A)$ where $\alpha$ is a
complex unit number. Our work generalizes some results obtained in
\cite{anradha2013, cuihou13, shader}.
\end{abstract}


\end{frontmatter}

\section{Introduction}

Throughout this paper, all matrices are real or complex. The identity matrix
of order $n$ is denoted by $I_{n}$ and the transpose of a matrix $A$ by
$A^{T}$. Let $\Sigma$ be a set of nonzero real or complex numbers such that if
$x\in \Sigma$ then $x^{-1}\in \Sigma$.\ Two square matrices $A$ and $B$ are
$\Sigma$-\emph{diagonally similar} if $B=\Lambda^{-1}A\Lambda$ for some
nonsingular diagonal matrix $\Lambda$ with diagonal entries in $\Sigma$. A square
matrix $A$ is \emph{reducible }if there exists a permutation matrix $P$, so
that $A$ can be reduced to the form $PAP^{T}=\left(
\begin{array}
[c]{cc}%
X & Y\\
0 & Z
\end{array}
\right)  $ where $X$ and $Z$ are square matrices. A matrix which is not
reducible is said to be \emph{irreducible}. A real matrix $A$ is nonnegative,
$A\geq0$, if all its entries are nonnegative.

Let $A$ be a $n\times n$ real or complex matrix. The multiset $\left \{
\lambda_{1},\lambda_{2},\ldots,\lambda_{n}\right \}  $ of eigenvalues of $A$ is
called the \emph{spectrum} of $A$ and is denoted by $sp(A)$. We usually assume
that $\left \vert \lambda_{1}\right \vert \geq \left \vert \lambda_{2}\right \vert
\geq \ldots \geq \left \vert \lambda_{n}\right \vert $. The \emph{spectral radius}
of $A$ is $\left \vert \lambda_{1}\right \vert $ that is denoted by $\rho(A)$.
Clearly, two real matrices $\left \{  -1,1\right \}  $-diagonally similar have
the same absolute value and the same spectrum.

In this paper, we address the following problem.

\begin{prob}
Let $A$ be a nonnegative real matrix and let $\alpha$ be a complex unit
number. Characterize the set of all signings of $A$ such that $sp(B)=\alpha sp(A)$.
\end{prob}

Our motivation comes from the works of Shader and So \cite{shader} which
correspond to adjacency matrix of a\emph{ }graph. Let $G$ be a finite simple
graph with vertex set $V(G)=\left \{  v_{1},\ldots,v_{n}\right \}  $ and edge
set $E(G)$. The \emph{adjacency matrix} of\emph{ }$G$ is the symmetric matrix
$A(G)=(a_{ij})_{1\leq i,j\leq n}$ where $a_{ij}=a_{ji}=1$\ if $\left \{
v_{i},v_{j}\right \}  $ is an edge of $G$, otherwise $a_{ij}=a_{ji}=0$. Since the matrix
$A(G)$ is symmetric, its eigenvalues  are real. The \emph{adjacency spectrum} $Sp\left(
G\right)  $ of $G$ is defined as the spectrum of $A(G)$. Let $G^{\sigma}$ be
an orientation of $G$, which assigns to each edge a direction so that the resultant graph $G^{\sigma}$ becomes an oriented graph. The
\emph{skew-adjacency} matrix of $G^{\sigma}$ is the real skew symmetric matrix
$S(G^{\sigma})=(a_{ij}^{\sigma})_{1\leq i,j\leq n}$\ where $a_{ij}^{\sigma
}=-a_{ji}^{\sigma}=1$ if $(i,j)$ is an arc of $G^{\sigma}$ and $a_{ij}%
^{\sigma}=0$, otherwise. The \emph{skew-spectrum} $Sp\left(  G^{\sigma
}\right)  $ of $G^{\sigma}$ is defined as the spectrum of $S(G^{\sigma})$.
Note that $Sp(G^{\sigma})$ consists of only purely imaginary eigenvalues
because $S(G^{\sigma})$ is a real skew symmetric matrix.

There are several recent works about the relationship between $Sp(G^{\sigma})$
and $Sp(G)$ (see for example \cite{anradha2013, anuradha2014, cavers,
cuihou13, shader}). In the last paper, Shader and So  have obtained a result which is closely related to our work. To state this result, let $G$ be a bipartite
graph with bipartition $[ I,J]$. The orientation $\varepsilon_{I,J}$
that assigns to each edge of $G$ a direction from $I$ to $J$ is called the
\emph{canonical orientation}. Shader and So showed \cite{shader} that
$Sp\left(  G^{\varepsilon_{I,J}}\right)  =iSp\left(  G\right)  $.        Moreover,
they proved that a graph $G$ is bipartite if and only if $Sp\left(  G^{\sigma}\right)  =isp\left(  G\right)  $ for some orientation $G^{\sigma}$ of $G$.
For connected graphs, this result can be viewed as a particular case of
Proposition \ref{shader generalis} because a graph $G$ is bipartite if and
only if there exists a permutation matrix $P$ such that $PA(G)P^{T}=\left(
\begin{array}
[c]{cc}%
0 & X\\
X^{T} & 0
\end{array}
\right)  $ where the zero diagonal blocks are square. Consider now an orientation
$G^{\sigma}$ of a graph $G$ and let $W$ be a subset of $V(G)$. The orientation
$G^{\tau}$ of $G$ obtained from $G^{\sigma}$ by reversing the direction of
all arcs between $W$ and $V(G)\backslash W$ is said to be obtained from
$G^{\sigma}$ by  \emph{switching} with respect to $W$. Two orientations are
\emph{switching-equivalent} if one can be obtained from the other by
switching. Anuradha, Balakrishnan, Chen, Li, Lian and So \cite{anradha2013}
proved that $Sp\left(  G^{\sigma}\right)  =isp\left(  G\right)  $ if an only
if $G^{\sigma}$ is switching-equivalent to the canonical orientation. For a
bipartite connected graph, this is a direct consequence of Proposition
\ref{solution shader} because two orientations of a graph are
switching-equivalent if and only if their skew-adjacency matrices are
$\{-1,1\}$-diagonally similar.

In this work, we consider only the case of irreducible matrices (not
necessarily symmetric). To state our main result, we need some terminology. A
\emph{digraph }$D$ is a pair consisting of a finite set $V(D)$ of
\emph{vertices} and a subset $E(D)$ of ordered pairs of vertices called
\emph{arcs}. Let $v,v^{\prime}$ two vertices of $D$, a \emph{path} $P$ from
$v$ to $v^{\prime}$ is a finite sequence $v_{0}=v,\ldots,v_{k}=v^{\prime}$
such that $(v_{0},v_{1}),\ldots,(v_{k-1},v_{k})$ are arcs of $D$. The
\emph{length} of $P$ is the number $k$ of its arcs.\ If $v_{0}=v_{k}$, we say
that $P$ is a \emph{closed path. }A digraph is said to be\emph{
strongly connected} if for any two vertices $v$ and $v^{\prime}$, there is a
path joining $v$ to $v^{\prime}$. It is easy to see that a strongly connected
digraph has a closed path. The \emph{period} of a digraph is the greatest
common divisor of the lengths of its closed paths. A digraph is
\emph{aperiodic }if its period is one.

With each $n\times n$ matrix $A=(a_{ij})_{1\leq i,j\leq n}$, we associate a
digraph $D_{A}$ on the vertex set $\left[  n\right]  =\left \{  1,\ldots
,n\right \}  $ and with arc set $E(D_{A})=\left \{  (i,j):a_{ij}\neq0\right \}
$. It is easy to show that $A$ is irreducible if and only if $D_{A}$ is
strongly connected. The \emph{period} of a matrix is the period of its
associate digraph.

\begin{nt}
Let $A$ be an irreducible nonnegative real matrix with period $p$ and let
$\alpha$ be a complex unit number, we denote by $\mathcal{M}(\alpha,A)$ the
set of all signings of $A$ such that $sp(B)=\alpha sp(A)$.
\end{nt}

In Corollary \ref{forme de alpha}, we gave a necessary condition for
$\mathcal{M}(\alpha,A)$ to be nonempty. More precisely, we will prove that if
$\mathcal{M}(\alpha,A)$ is nonempty, then
$\alpha=e^{\frac{i\pi k}{p}}$ for some $k\in \left \{  0,\ldots,2p-1\right \}  $.
If $p>1$, we will construct (see Corollary \ref{nonempty period p}) for each
$k\in \left \{  0,\ldots,2p-1\right \}  $ a matrix $\widetilde{A}$ such that
$sp(\widetilde{A})=e^{\frac{i\pi k}{p}}sp(A)$. Finally, by using Corollary
\ref{forme de alpha} and Proposition \ref{structurede malpha}, we prove our main theorem.

\begin{thm}
\label{main theorem}Let $A$ be an irreducible nonnegative matrix with period
$p$ and let $\alpha$ be an unit complex number. Then the following statements hold:

\begin{description}
\item[i)] $\mathcal{M}(\alpha,A)$ is nonempty iff $\alpha=e^{\frac{i\pi k}{p}%
}$ for some $k\in \left \{  0,\ldots,2p-1\right \}  $;

\item[ii)] For $k\in \left \{  0,\ldots,2p-1\right \}  $, we have

\item[a)] if $k$ is even then $\mathcal{M}(e^{\frac{i\pi k}{p}},A)$ is the set
of all matrices $\left \{  -1,1\right \}  -$diagonally similar to $A$;

\item[b)] if $k$ is odd then $\mathcal{M}(e^{\frac{i\pi k}{p}},A)$ is the set
of all matrices $\left \{  -1,1\right \}  -$diagonally similar to $B_{0}$ where
$B_{0}$ is an arbitrary  signing of $A$ such that  $sp(B_{0})=e^{\frac{i\pi}{p}}sp(A)$.
\end{description}
\end{thm}

\begin{req}
Let $A$ be an aperiodic matrix. It follows from the above Theorem
that $\mathcal{M}(\alpha,A)$ is empty whenever $\alpha \neq \pm1$. Moreover,
$\mathcal{M}(1,A)$ (resp. $\mathcal{M}(-1,A)$) is the set of all matrices
$\left \{  -1,1\right \}  -$diagonally similar to $A$ (resp. to $-A$).
\end{req}

\section{Some properties of $\mathcal{M}(\alpha,A)$}

We will use the following theorem due to Frobenius (see \cite{godsil}).

\begin{thm}
\label{frobenius} Let $A$ be an irreducible nonnegative $n\times n$ matrix and
let $B$ be a complex $n\times n$ matrix such that $|B|\leq A$. Then
$\left \vert \lambda \right \vert \leq \rho(A)$ for each eigenvalue $\lambda$ of
$B$ and $\lambda=\rho(A)e^{i\theta}$ iff $B=e^{i\theta}LAL^{-1}$, where $L$ is
a complex diagonal matrix such that  $\left \vert L\right \vert =I_{n}$.
\end{thm}

In the next Proposition, we describe the structure of $\mathcal{M}(\alpha,A)$.

\begin{prop}
\label{structurede malpha} Let $A$ be an irreducible nonnegative matrix with
period $p$ and let $\alpha$ be an unit complex number. Then $\mathcal{M}%
(\alpha,A)$ is empty or $\mathcal{M}(\alpha,A)$ is the set of all  matrices
$\left \{  -1,1\right \}  -$diagonally similar to $B_{0}$ where $B_{0}$ is an
arbitrary matrix in $\mathcal{M}(\alpha,A)$.
\end{prop}

\begin{pf}
Assume that $\mathcal{M}(\alpha,A)$ is nonempty and let $B_{0}\in
\mathcal{M}(\alpha,A)$. It is easy to see that $\Lambda^{-1}B_{0}\Lambda\in \mathcal{M} (\alpha,A)$ for every $\left \{  -1,1\right \}  $-diagonal matrix $\Lambda$.
Conversely, let $B\in \mathcal{M}(\alpha,A)$. Then $sp(B)=sp(B_{0})=\alpha
sp(A)$. Let $\lambda=\rho(A)e^{i\theta}$ be a common eigenvalue of $B$ and
$B_{0}$. By Theorem \ref{frobenius}, we have $B=e^{i\theta}LAL^{-1}$ and
$B_{0}=e^{i\theta}L^{\prime}AL^{\prime-1}$ where $L$, $L^{\prime}$ are complex
diagonal matrices such that $\left \vert L\right \vert =\left \vert L^{\prime
}\right \vert =I_{n}$. It follows that $B=(L^{^{\prime}}L)^{-1}B_{0}L^{\prime
}L^{-1}$. To conclude, it suffices to apply Lemma \ref{realsimilarit} below.
\qed\end{pf}

\begin{lem}
\label{realsimilarit} Let $A$  be an  $n\times n$  irreducible
nonnegative matrix and let $B$, $B^{\prime}$  be two signings of $A$. If there exists
a complex diagonal matrix $\Gamma$ such that $B^{\prime}=\Gamma B\Gamma^{-1}$
and $\left \vert \Gamma \right \vert =I_{n}$ then $B$ and $B^{\prime}$ are
$\left \{  -1,1\right \}  -$diagonally similar.
\end{lem}

\begin{pf}
Let $A:=(a_{ij})_{1\leq i,j\leq n}$, $B:=(b_{ij})_{1\leq i,j\leq n}$ and  $B^{\prime}:
=(b_{ij}^{\prime})_{1\leq i,j\leq n}$. We denote by $\gamma_{1},\ldots,\gamma_{n}$ the diagonal entries of $\Gamma$. Let $\Delta$ be the diagonal matrix with diagonal entries $\delta_{j}%
=\gamma_{j}\gamma_{1}^{-1}$ for $j=1,\ldots, n$. As $A$ is irreducible, the
digraph $D_{A}$ is strongly connected and then there is a path $j=i_{1}%
,\ldots,i_{r}=1$ of $D_{A}$ from $j$ to $1$. By definition of $D_{A}$, we have
$a_{i_{1}i_{2}}\neq0,\ldots,a_{i_{r-1}i_{r}}\neq0$. It follows that
$b_{i_{1}i_{2}}\neq0$,\ldots, $b_{i_{r-1}i_{r}}\neq0$. and $b_{i_{1}i_{2}%
}^{\prime}\neq0$,\ldots, $b_{i_{r-1}i_{r}}^{\prime}\neq0$ because $\left \vert
B\right \vert =\left \vert B^{\prime}\right \vert =A$. Moreover, from the
equality $B^{\prime}=\Gamma B\Gamma^{-1}$ we have $b_{i_{1}i_{2}}^{\prime
}=\gamma_{i_{1}}b_{i_{1}i_{2}}\gamma_{i_{2}}^{-1}$ ,$b_{i_{2}i_{3}}^{\prime
}=\gamma_{i_{2}}b_{i_{2}i_{3}}\gamma_{i_{3}}^{-1}$,\ldots, $b_{i_{r-1}i_{r}%
}^{\prime}=\gamma_{i_{r-1}}b_{i_{r-1}i_{r}}\gamma_{i_{r}}^{-1}$. Then
$b_{i_{1}i_{2}}^{\prime}\ldots b_{i_{r-1}i_{r}}^{\prime}=\gamma_{i_{1}}%
\gamma_{i_{r}}^{-1}b_{i_{1}i_{2}}\ldots b_{i_{r-1}i_{r}}$. But by hypothesis,
$B,B^{\prime}$ are real matrices and $\left \vert B\right \vert =\left \vert
B^{\prime}\right \vert $, then $b_{i_{1}i_{2}}^{\prime}\ldots b_{i_{r-1}i_{r}%
}^{\prime}=\pm b_{i_{1}i_{2}}\ldots b_{i_{r-1}i_{r}}$ and hence $\delta
_{j}=\gamma_{j}\gamma_{1}^{-1}=\gamma_{i_{1}}\gamma_{i_{r}}^{-1}\in \left \{
-1,1\right \}  $. To conclude, it suffices to see that $\Delta B\Delta
^{-1}=\Gamma B\Gamma^{-1}=B^{\prime}$.
\qed\end{pf}

By using Theorem \ref{frobenius}, we obtain the following.

\begin{prop}
\label{eigenvalue of B} Let $A$ be an irreducible
nonnegative real matrix with period $p$ and let $B$
be a signing of $A$ such that $\rho(B)=\rho(A)$. If $\lambda$ is an eigenvalue
of $B$ such that $\left \vert \lambda \right \vert =\rho(A)$,  then $\lambda
=\rho(A)e^{\frac{i\pi k}{p}}$ for some $k\in \left \{  0,\ldots,2p-1\right \}  $.
\end{prop}

\begin{pf}
 Let $A:=(a_{ij})_{1\leq i,j\leq n}$, $B:=(b_{ij})_{1\leq i,j\leq n}$ and
 $\lambda=\rho(A)e^{i\theta}$. By Theorem \ref{frobenius}, we have
$B=e^{i\theta}LAL^{-1}$ where $L$ is a complex diagonal matrix such that  $\left \vert
L\right \vert =I_{n}$. It follows that $b_{ij}=e^{i\theta}l_{i}a_{ij}l_{j}%
^{-1}$ for $i,j\in \left \{  1,\ldots,n\right \}  $, where $l_{1},\ldots,l_{n}$
are the diagonal entries of $L$. Consider now a closed path $C=(i_{1}%
,i_{2},\ldots,i_{r},i_{1})$ of $D_{A}$. By the previous equality, we have%
\begin{align*}
\text{ }\frac{b_{i_{1}i_{2}}\ldots b_{i_{r-1}i_{r}}b_{i_{r}i_{1}}}%
{a_{i_{1}i_{2}}\ldots a_{i_{r-1}i_{r}}a_{i_{r}i_{1}}}  &  =(e^{i\theta
}l_{i_{1}}l_{i_{2}}^{-1})\ldots(e^{i\theta}l_{i_{r-1}}l_{i_{r}}^{-1}%
)(e^{i\theta}l_{i_{r}}l_{i_{1}}^{-1})\\
&  =(e^{i\theta})^{r}%
\end{align*}

Then $(e^{i\theta})^{r}\in \left \{  1,-1\right \}  $ because $\left \vert
B\right \vert =A$.

Since $p$ is the greatest common divisor of the lengths of the closed paths in
$D_{A}$, we have $(e^{i\theta})^{p}\in \left \{  1,-1\right \}  $ and then
$\lambda=\rho(A)e^{\frac{i\pi k}{p}}$ for some $k\in \left \{  0,\ldots
,2p-1\right \}  $.
\qed\end{pf}

\begin{req}
\label{eigenvalues of A} Let $A=(a_{ij})_{1\leq i,j\leq n}$ be an irreducible
nonnegative real matrix with period $p$ and let $\lambda$ be an eigenvalue of
$A$ such that $\left \vert \lambda \right \vert =\rho(A)$. By applying
Proposition \ref{eigenvalue of B} to $B=A$, we have $\lambda=\rho
(A)e^{\frac{i\pi k}{p}}$ for some $k\in \left \{  0,\ldots,2p-1\right \}  $.
\end{req}

The following result gives a necessary condition under which $\mathcal{M}%
(\alpha,A)$ is nonempty.

\begin{cor}
\label{forme de alpha} Let $A$ be an irreducible nonnegative real matrix of
period $p$ and let $\alpha$ be a complex unit number. If $\mathcal{M}%
(\alpha,A)$ is nonempty then $\alpha=e^{\frac{i\pi k}{p}}$ for some
$k\in \left \{  0,\ldots,2p-1\right \}  $, in particular $\alpha^{p}=\pm1$.
\end{cor}

\begin{pf}
Let $\lambda$ be an eigenvalue of $A$ such that $\left \vert \lambda \right \vert
=\rho(A)$. By applying Remark \ref{eigenvalues of A}, we have $\lambda
=\rho(A)e^{\frac{i\pi k}{p}}$ for some $k\in \left \{  0,\ldots,2p-1\right \}  $.
Let $B\in \mathcal{M}(\alpha,A)$. We have $\alpha \rho(A)e^{\frac{i\pi k}{p}}\in
sp(B)$ because $sp(B)=\alpha sp(A)$. It follows from Proposition
\ref{eigenvalue of B} that $\alpha \rho(A)e^{\frac{i\pi k}{p}}=\rho
(A)e^{\frac{i\pi h}{p}}$ for some $h\in \left \{  0,\ldots,2p-1\right \}  $ and
hence $\alpha=e^{\frac{i\pi(h-k)}{p}}$.
\qed\end{pf}

\section{Proof of the Main Theorem}

Let $n$ be a positive integer and let $(r_{1},\ldots,r_{p})$ be a partition of
$n$, that is $r_{1},\ldots,r_{p}$ are positive integers and $r_{1}%
+\cdots+r_{p}=n$. For $i=1,\ldots,p-1$, let $A_{i}$ be a $r_{i}\times r_{i+1}%
$ matrix and let $A_{p}$ be a $r_{p}\times r_{1}$ matrix. The matrix
$\left(
\begin{array}
[c]{ccccc}%
0 & A_{1} & 0 & \cdots & 0\\
0 & 0 & A_{2} & \cdots & 0\\
\vdots & \vdots & \ddots & \ddots & \vdots \\
0 & 0 & 0 & \ddots & A_{p-1}\\
A_{p} & 0 & \cdots & 0 & 0
\end{array}
\right)  $ is denoted by $Cyc(A_{1},A_{2},\ldots,A_{p})$.

Each matrix of this form is called $p$-\emph{cyclic. }

Recall\ the well-known result of Frobenius about irreducible matrices with
period $p>1$.

\begin{prop}
\label{cyclic form} Let $A$ be an irreducible nonnegative real matrix with
period $p>1$, then there exists a permutation matrix $P$ such that $PAP^{T}$
is $p$-cyclic.
\end{prop}

\begin{req}
\label{coverce aperiodic} For aperiodic matrices, the converse of Corollary
\ref{forme de alpha} is true because $A\in \mathcal{M}(1,A)$ and $-A\in
\mathcal{M}(-1,A)$.
\end{req}

For $p>1$, we have the following result.

\begin{prop}
\label{converse for periodic} Let $A=Cyc(A_{1},A_{2},\ldots,A_{p})$ be a
nonnegative $p$-cyclic matrix where $A_{i}$ is a $r_{i}\times r_{i+1}$ matrix
for $i=1,\ldots,p-1$ and $A_{p}$ is a $r_{p}\times r_{1}$ matrix. Let
$\widetilde{A}$ be the matrix obtained from $A$ by replacing the block $A_{p}$
by $-A_{p}$. Given $\alpha=e^{\frac{i\pi k}{p}}$ where $k\in \left \{
0,\ldots,2p-1\right \}  $, then

\begin{description}
\item[i)] if $k$ is even, $e^{\frac{i\pi k}{p}}A$ is diagonally similar to
$A$, in particular $sp(A)=e^{\frac{i\pi k}{p}}sp(A)$;

\item[ii)] if $k$ is odd, $e^{\frac{i\pi k}{p}}A$ is diagonally similar to
$\widetilde{A}$, in particular $sp(\widetilde{A})=e^{\frac{i\pi k}{p}}sp(A)$.
\end{description}
\end{prop}

\begin{pf}
Let $L:=\left(
\begin{array}
[c]{ccccc}%
I_{r_{1}} & 0 & 0 & \cdots & 0\\
0 & e^{\frac{i\pi}{p}}I_{r_{2}} & 0 & \cdots & 0\\
\vdots & \vdots & \ddots & \ddots & \vdots \\
0 & 0 & 0 & \ddots & 0\\
0 & 0 & \cdots & 0 & e^{\frac{i\pi(p-1)}{p}}I_{r_{p}}%
\end{array}
\right)  $. It easy to check that if $k$ is even, $e^{\frac{i\pi k}{p}%
}LAL^{-1}=A$ and if $k$ is odd, $e^{\frac{i\pi k}{p}}LAL^{-1}=\widetilde{A}$.
\qed\end{pf}

The next Corollary is a direct consequence of the above Proposition and
Theorem \ref{cyclic form}.

\begin{cor}
\label{nonempty period p} Let $A$ be an irreducible nonnegative matrix with
period $p>1$. Then $\mathcal{M}(e^{\frac{i\pi k}{p}},A)$ is nonempty for
$k\in \left \{  0,\ldots,2p-1\right \}  $.
\end{cor}

\begin{pf}
By Proposition \ref{cyclic form}, there exists a permutation matrix $P$ such
that $PAP^{T}$ is $p-$cyclic. Let $A^{\prime}:=PAP^{T}:=Cyc(A_{1},A_{2}%
,\ldots,A_{p})$ and let $\widetilde{A^{\prime}}$ the matrix obtained from
$A^{\prime}$ by replacing the block $A_{p}$ by $-A_{p}$. It follows from
Proposition \ref{converse for periodic} that $sp(A^{\prime})=e^{\frac{i\pi
k}{p}}sp(A^{\prime})$ if $k$ is even and $sp(\widetilde{A^{\prime}}%
)=e^{\frac{i\pi k}{p}}sp(A^{\prime})$ if $k$ is odd. If $k$ is even then
$sp(A)=e^{\frac{i\pi k}{p}}sp(A)$ because $P^{T}A^{\prime}P=A$ and hence
$A\in \mathcal{M}(e^{\frac{i\pi k}{p}},A)$. If $k$ is odd then $sp(P^{T}%
\widetilde{A^{\prime}}P)=e^{\frac{i\pi k}{p}}sp(P^{T}A^{\prime}P)=e^{\frac
{i\pi k}{p}}sp(A)$. Moreover, since$\left \vert \widetilde{A^{\prime}%
}\right \vert =A^{\prime}$, we have $\left \vert P^{T}\widetilde{A^{\prime}%
}P\right \vert =P^{T}A^{\prime}P=A$ and then $P^{T}\widetilde{A^{\prime}}%
P\in \mathcal{M}(e^{\frac{i\pi k}{p}},A)$.
\qed\end{pf}

\begin{req}
\label{egalite des malpha} It follows from Corollary \ref{forme de alpha},
Remark \ref{coverce aperiodic} and Corollary \ref{nonempty period p} that
$\mathcal{M}(\alpha,A)$ is nonempty iff $\alpha=e^{\frac{i\pi k}{p}}$ for some
$k\in \left \{  0,\ldots,2p-1\right \}  $. Moreover, as $e^{\frac{i\pi k}{p}%
}sp(A)=sp(A)$ if $k$ is even, we have $e^{\frac{i\pi k}{p}}sp(A)=e^{\frac
{i\pi}{p}}sp(A)$ if $k$ is odd and then
\begin{align*}
\mathcal{M}(1,A)  &  =\mathcal{M}(e^{\frac{2i\pi}{p}},A)=\cdots=\mathcal{M}%
(e^{\frac{2(p-1)i\pi}{p}},A)\\
\mathcal{M}(e^{\frac{i\pi}{p}},A)  &  =\mathcal{M}(e^{\frac{3i\pi}{p}%
},A)=\cdots=\mathcal{M}(e^{\frac{(2p-1)i\pi}{p}},A)
\end{align*}

\end{req}

\begin{pot}

i) (Necessity) This follows from Corollary \ref{forme de alpha}.

(Sufficiency) If $p=1$ then from Corollary \ref{forme de alpha} we have
$\alpha=\pm1$ and hence by Remark \ref{coverce aperiodic}, $\mathcal{M}%
(\alpha,A)$ is nonempty. If $p>1$, it suffices to apply Corollary
\ref{nonempty period p}.

ii) By Remark \ref{egalite des malpha}, $\mathcal{M}(e^{\frac{i\pi k}{p}%
},A)=\mathcal{M}(1,A)$ if $k$ is even and $\mathcal{M}(e^{\frac{i\pi k}{p}%
},A)=\mathcal{M}(e^{\frac{i\pi}{p}},A)$ if $k$ is odd. To prove a) and b) it
suffices to apply Proposition \ref{structurede malpha} respectively to
$\alpha=1$ and to $\alpha=e^{\frac{i\pi}{p}}$.
\qed\end{pot}

\section{The special case of symmetric matrices}

In this section we give some consequences of Theorem \ref{main theorem} when
the matrix $A$ is symmetric.

\begin{prop}
\label{shader generalis} Let $A$ be an irreducible nonnegative symmetric real
matrix. Then the following statements are equivalent :

\begin{description}
\item[i)] there is a real skew-symmetric matrix $B$ such that $\left \vert
B\right \vert =A$ and $sp(B)=isp(A)$;

\item[ii)] there exists a permutation matrix $P$ such that $PAP^{T}=\left(
\begin{array}
[c]{cc}%
0 & X\\
X^{T} & 0
\end{array}
\right)  $ where the zero diagonal blocks are square.
\end{description}
\end{prop}

\begin{pf}
Note that the period of  $A$ is a most $2$ because its associate digraph
contains a closed path of length $2$.

First, we prove the implication i)$\implies$ii). It follows from i) that
$\mathcal{M}(i,A)$ is nonempty. Then, by Corollary \ref{forme de alpha}, the
period of $A$ is necessarily $2$. It follows from Proposition
\ref{cyclic form} that there exists a permutation matrix $P$, so that
$PAP^{T}=\left(
\begin{array}
[c]{cc}%
0 & X\\
Y & 0
\end{array}
\right)  $ where the zero diagonal blocks are square. But, the matrix
$PAP^{T}$ is symmetric, then $Y=X^{T}$.

To prove ii) $\implies$ i), it suffices to apply Proposition
\ref{converse for periodic} for $p=2$ and $k=1$.
\qed\end{pf}

\begin{prop}
\label{solution shader} Let $A=\left(
\begin{array}
[c]{cc}%
0 & X\\
X^{T} & 0
\end{array}
\right)  $ be an irreducible nonnegative symmetric matrix and let $B$ be a
skew symmetric matrix such that $\left \vert B\right \vert =A$. Then the
following statements are equivalent :

\begin{description}
\item[i)] $sp(B)=isp(A)$;

\item[ii)] $B$ is $\{-1,1\}$-diagonally similar to $\widetilde{A}=\left(
\begin{array}
[c]{cc}%
0 & X\\
-X^{T} & 0
\end{array}
\right)  $.
\end{description}
\end{prop}

\begin{pf}
It follows from assertion ii) of Proposition \ref{converse for periodic} that
$sp(\widetilde{A})=isp(A)$ and then $\widetilde{A}\in \mathcal{M}(i,A)$. The
equivalence i) $\Longleftrightarrow$ ii) result from assertion b) of Theorem
\ref{main theorem}.
\qed\end{pf}

\end{document}